\documentclass[11pt]{amsart}
\usepackage[dvipsnames]{xcolor}
\newcommand{\ben}{\begin{enumerate}}
\newcommand{\een}{\end{enumerate}}
\newcommand{\eq}[2][label]{\begin{equation}\label{#1}#2\end{equation}}
\newcommand{\av}[2]{\langle #1\rangle_{_{\scriptstyle #2}}}

\newcommand{\ve}{\varepsilon}

\usepackage{amsfonts}
\usepackage{amssymb}
\usepackage{amsmath}
\usepackage{graphicx,enumitem,url,hyperref}

\newcommand{\BMO}{{\rm BMO}}
\newcommand{\VMO}{{\rm VMO}}

\newcommand{\Cut}{{\rm Cut}}

\newcommand{\Xicl}[1]{{\mathcal L}_{#1}}

\renewcommand{\phi}{\varphi}

\newcommand{\AAA}{A}

\newtheorem{theorem}{Theorem}[section]

\newtheorem{lemma}[theorem]{Lemma}
\newtheorem{corollary}[theorem]{Corollary}
\newtheorem{prop}[theorem]{Proposition}

\newtheorem*{theorem*}{Theorem}{\bf}{\it}
\newtheorem*{proposition*}{Proposition}{\bf}{\it}
\newtheorem*{observation*}{Observation}{\bf}{\it}
\newtheorem*{lemma*}{Lemma}{\bf}{\it}
\newtheorem*{conjecture*}{Conjecture}{\bf}{\it}

\theoremstyle{definition}
\newtheorem{definition}[theorem]{Definition}

\theoremstyle{remark}
\newtheorem{remark}[theorem]{Remark}

\newcommand{\kk}[1]{k\left(#1\right)}
\newcommand{\rr}[1]{r\left(#1\right)}

\newcommand{\kr}[1]{\frac{k(#1)}{r(#1)}}
\newcommand{\rk}[1]{\frac{r(#1)}{k(#1)}}

\newcommand{\Cnn}{\Gamma}
\newcommand{\cnn}{\gamma}
\newcommand{\de}[1]{\left[#1\right]'}
\newcommand{\Cn}[3]{\Gamma_{#3}^{#2}(#1)}
\newcommand{\cn}[3]{\gamma_{#3}^{#2}(#1)}
\newcommand{\F}[2]{F_2^{#1}(#2)}
\newcommand{\AP}[2]{\av{#1}{#2}\av{#1^{-1}}{#2}}

\newcommand{\UU}[2]{U(#1_1,#1_2,#2)}
\newcommand{\TT}[2]{\mathcal{T}(#1_1,#1_2,#2)}
\newcommand{\M}{M}

\DeclareMathOperator{\col}{Cut}

\numberwithin{equation}{section}

\allowdisplaybreaks
\begin{document}
\title
{Monotone rearrangement for flattening weights $F_2$}

\author{Sasanka Wijesiri Gunawardana}
\date{May 1, 2026}
\begin{abstract}
We consider a subclass of the Muckenhoupt $A_2$ weights with added smoothness, which we refer to as flattening weights $F_2.$ We quantify how such weights arise as exponents of VMO functions and prove that monotone rearrangement does not increase the class modulus~in~$F_2.$   
\end{abstract}
\maketitle
\section{Introduction and main results.}
\label{I}
\subsection{Preliminaries}
Throughout this paper, $I$ and $J$ represent finite intervals; the symbol $\av{\cdot }{J}$ stands for the average over $J$ with respect to the Lebesgue measure.  A weight is an almost everywhere positive function. A key class of weights is the $A_p$, introduced by Muckenhoupt in~\cite{muck}. Here we need the definition of $A_2$ restricted to an interval $I$:
$$
w\in A_2(I)\quad\Longleftrightarrow\quad [w]_{A_2(I)}:=\sup_{J\subset I}\av{w}J\av{w^{-1}}J<\infty.
$$
The number $[w]_{A_2(I)}$ is called the $A_2$-characteristic of $w$ on $I.$ Note that $[w]_{A_2(I)}=1$ if and only if $w$ is constant a.e. on $I.$

It is classical that $A_2$ weights (and in general $A_p$ weights for $p>1$) are exponents of $\BMO$ functions (functions of bounded mean oscillation). Recall that a function $\varphi$ on $I$ is said to be in $\BMO(I)$ if and only if
$$
\|\varphi\|_{\BMO(I)}:=\sup_{J\subset I}\big(\av{\varphi^2}J-\av{\varphi}J^2\big)^{1/2}<0.
$$ 
The following result can be found in, e.g., \cite{duo}:
\begin{prop}
\label{prop1}
~
\ben
\item
If $w\in A_2(I),$ then $\log w\in\BMO.$
\item
Conversely, if $\varphi\in\BMO(I)$ and $\lambda>0$ is small enough, then $e^{\lambda \varphi}\in A_2(I).$
\een
\end{prop}
An important subspace of $\BMO$ is VMO, the space of functions of vanishing mean oscillation. It was introduced by Sarason in~\cite{sarason}. A function $\varphi\in\BMO(I)$ is said to be in $\VMO(I)$ if
$$
\sup_{J\subset I, |J|\le\delta}\big(\av{\varphi^2}J-\av{\varphi}J^2\big)\to 0,~\text{as }\delta\to0.
$$
Sarason also proved the following theorem:
\begin{theorem}[\cite{sarason}]
Let $\varphi\in\BMO(I).$ Then $\varphi\in\VMO(I)$ if and only if
\eq[sar1]{
\sup_{J\subset I, |J|\le\delta} \av{e^\varphi}J\av{e^{-\varphi}}J\to 1,~\text{as }\delta\to0.
}
\end{theorem}
One can view condition~\eqref{sar1} as a requirement that the weight $e^\varphi$ become closer and closer to a constant (i.e., ``flatten'') as the size of the interval over which its $A_2$ characteristic is computed goes to 0. We are interested in quantitating this notion and in further analysis of such weights. To this end, we need a quantitative version of VMO; the following definition comes from \cite{sz1}.  

Fix an interval $I$ and let $\xi \colon [0,|I|] \to [0,\infty)$ be an increasing continuous function such that
$\xi(0)=0.$ We will call such functions $\VMO$-moduli. Let $\Xicl{\xi}(I)$ be the class of all functions $\varphi\colon I\to\mathbb{R}$ for which the following quantity is finite:
\eq[norm]{
\|\phi\|_{\Xicl{\xi}(I)} := \inf\big\{C>0 \colon \av{\phi^2}{J} - \av{\phi}{J}^2 \leq C^2 \xi^2(|J|), J \subset I\big\}.
}
To define a flattening version of $A_2,$ let $k\colon [0,\infty)\to[1,|I|]$ be an increasing continuous function such that $\kk{0}=1$. Let us call such functions $F_2$-moduli.  
\begin{definition}
For an interval $I$ and an $F_2$-modulus $k$ on $I,$ we say that a weight $w$ belongs to the class $\F{k}{I}$ if
\begin{equation}
\label{f2}
    \AP{w}{J}\le k(|J|)\text{\space for all\space} J\subset I.
\end{equation}
\end{definition}
Clearly, $\F{k}I\subset A_2(I);$ indeed, if $w\in\F{k}I,$ then $[w]_{A_2(I)}\le k(|I|).$ 
\begin{remark}
In what follows we will simply write ``modulus" for $F_2$-modulus, unless we need to distinguish it from a $\VMO$-modulus.
\end{remark}

We now state an analog of the $A_2$-$\BMO$ relationship from Proposition~\ref{prop1} for the classes $F^k_2$ and $\Xicl{\xi},$ one that captures the sharp relationship between the $F_2$-modulus $k$ and the $\VMO$-modulus $\xi.$
\begin{theorem}~
\label{IT2}
\ben[leftmargin=.25in]
\item
    If $w\in \F{k}{I}$, then $\log w\in \Xicl{\sqrt{k-1}}(I)$ and $\|\log w\|_{\Xicl{\sqrt{k-1}}(I)}\le C$ for an absolute constant $C.$
    \item 
Conversely, if $\varphi\in \mathcal{L}_{\xi}(I),$ $\|\varphi\|_{\Xicl{\xi}}(I)=1,$
    and
    $\lambda\in[0,\frac1{\xi(|I|)})$, then $e^{\lambda \varphi}\in \F{1+\lambda^2M\xi^2}{I},$ for some constant $M$ depending only on the product $\lambda\xi(|I|).$
\een
\end{theorem}
\begin{remark}
The proof of Theorem~\ref{IT2}, given in Section~\ref{proofs}, yields explicit expressions for $M$ and $C.$ However, neither is sharp, so we do not include them in the statement.
\end{remark}

We also record a simple result concerning the content of the class $F_2^k.$ The following result was proved in~\cite{sz1}.
\begin{prop}[\cite{sz1}]
\label{oPL2}
The class $\Xicl{\xi}(I)$ contains non-constant functions if and only if
\begin{equation*}
    \liminf_{t\to0^+} \frac{\xi(t)}{t}>0.
\end{equation*}
    \end{prop}
The analogous statement for $F_2^k$ is a corollary of Proposition~\ref{oPL2} and Theorem~\ref{IT2}. Its proof is also given in Section~\ref{proofs}.
\begin{corollary}
\label{IT1}
    The class $\F{k}{I}$ contains non-constant functions if and only if
    \begin{equation}
    \label{NC}
        \liminf_{t\rightarrow 0^+}\, \frac{\sqrt{k(t)-1}}{t}>0.
    \end{equation}
\end{corollary}

\subsection{Main result}
We now turn to our main result, which states that the class $F_2^k(I)$ is invariant under monotone rearrangement. To define this operation, take an interval $I=(a,b)$ and a weight $w\colon I\to [0,\infty).$ The decreasing rearrangement of $w$ is the weight $w^*\colon I\to [0,\infty)$ given by
$$
w^*(t)=\inf\{\lambda:|\{s\in I\colon w(s)>\lambda\}|\le t-a\},\quad t\in I.
$$
As is well known, $w^*$ is the unique (up to sets of measure zero) decreasing function on $I$ that is equidistributed with $w$: 
$$
|\{s\in I:w^*(s)>\lambda\}|=|\{s\in I:w(s)>\lambda\}|,\quad \text{for all }\lambda\ge0.
$$
\begin{theorem}
\label{MT1} Take an interval $I$ and let $k$ be a modulus on $I.$ If $w\in \F{k}{I}$, then $w^*\in \F{k}{I}$.
\end{theorem}
\bigskip

The motivation behind Theorem~\ref{MT1} is two-fold. The first goal is to add another, fairly delicate result to the study of the behavior of averaging classes under rearrangement. This line of inquiry has a rich history in literature. Klemes~\cite{Klemes1} showed that rearrangement does not increase the $\BMO^1$ norm in dimension 1. A detailed investigation of various rearrangement techniques in $\BMO^1,$ as well as in weight classes on rectangles, can be found in the book~\cite{Kor2}. The article \cite{studia} gives the sharp constant of rearrangement from the dyadic $\BMO([0,1]^n)$ to the usual BMO on an interval. One can also study the continuity of the rearrangement operator and in this direction we mention the article~\cite{Ryan2}, where it is proved that the rearrangement is continuous as an operator from VMO to VMO. A key general result was obtained by Zatitskii and Stolyarov in~\cite{sz}, where they proved that the rearrangement preserves general classes $A_\Omega(I).$ This formulation includes BMO with the quadratic norm and $A_p,$ $1<p\le\infty$; their result means that the BMO norm and $A_p$ characteristic do not increase under rearrangement. 

Studying classes without scale invariance, such as VMO and $F_2,$ is naturally more difficult. Recently, Slavin and Zatitskii proved in~\cite{sz1} that monotone rearrangement does not increase the $\Xicl{\xi}$ norm. We adopt their methods in the present paper.

Our second motivation is the future development of a theory of sharp integral estimates for the class $F_2$ via Bellman functions. The application of that method in scale-dependent settings is very recent (see the preprint~\cite{sz2}) and it currently requires working with monotone functions, which in the $F_2$ case will be enabled by Theorem~\ref{MT1}.
\bigskip

The rest of the paper is devoted to the proof of Theorem~\ref{MT1}, except for Section~\ref{proofs}, where we prove Theorem~\ref{IT2} and Corollary~\ref{IT1}. We present the proof of our main theorem in three stages, relaxing requirements from one stage to the next.

\subsection*{Smooth convex case} 
\label{smm}
This case is given in Section~\ref{SmoSet} and takes the vast majority of the effort.
Here we prove Theorem~\ref{MT1} under additional assumptions. Specifically, let
\begin{equation}
\label{I12}
    \AAA(t)=t^2(k(t)-1).
\end{equation} 
In the smooth convex case, we assume that the modulus $k$ is twice continuously differentiable on $(0,|I|]$ and that
$
k'(t)>0, \AAA''(t)>0,
$
$t\in(0,|I|].$
\subsection*{General convex case} This case is presented in Section~\ref{gen1}. Here we drop the differentiability requirements on $k$ and $\AAA$ and simply assume that $\AAA$ is convex on $[0,|I|].$
\subsection*{General case} This case is presented in Section~\ref{gen2}. Here we prove Theorem~\ref{MT1} in full generality, as stated.

\section{Proof of the main  theorem in the smooth convex case}
\label{SmoSet}

Throughout this section, we are working under the assumptions stated above: $k$ is twice continuously differentiable on $(0,|I|]$, and $k'(t)>0,$ $\AAA''(t)>0$ for all $t\in(0,|I|];$ here $\AAA$ is given by~\eqref{I12}.

\subsection{The function $r$: definition and properties}
Define a function $r\colon [0,|I|]\to [1,\infty)$ by
\begin{equation}
  \label{r1}
  \rr{t}=\begin{cases}
     k+\frac12\,t k'(t)+\sqrt{\big(k+\frac12\,tk'(t)\big)^2-k(t)},& \text{if }t>0;\\
      1,& \text{if } t=0.
  \end{cases}
\end{equation}
From~\eqref{r1}, it is clear that for $t>0$ we have $r>1$ and that $r$ it is the larger solution of the quadratic equation
\begin{equation}
\label{r2}
    2+t\,\frac{k'(t)}{k(t)}=\frac{1}{\rr{t}}+\rk{t}.
\end{equation}
By replacing $k$ in \eqref{r2} by $r\,\frac{k}{r}$ and rearranging terms, we get
\begin{equation}
\label{r3}
    \frac{\de{t(\rr{t}-1)}}{\rr{t}}=-\rk{t}\bigg[t\bigg(\kr{t}-1\bigg)\bigg]'.
\end{equation}
\begin{lemma}
\label{rL2}
 The function $r$ is continuous on $[0,|I|].$
\end{lemma}
\begin{proof}
    From~\eqref{r1}, the continuity of $r$ for $t>0$ is obvious. To check the continuity at 0, observe that since the function $\AAA$ is convex, for $t>0$ we have
$$
A'(t)\le \frac{A(2t)-A(t)}t\le \frac{A(2t)}t=4t(k(2t)-1),
$$
or, equivalently,
$$
2(k(t)-1)+tk'(t)\le 4(k(2t)-1).
$$
Therefore, $\lim_{t\to0^+}tk'(t)=0,$ and from~\eqref{r1}, we have $\lim_{t\to0^+}r(t)=1.$
\end{proof}
\subsection{The special curves and their properties}
We first define the following domain in $\mathbb{R}^3$ that gives the class $F_2^k$ a geometric characterization.
\begin{equation}
\label{dm1}
    \omega=\left\{(x_1,x_2,t)\in\mathbb{R}^3\colon 0\le t\le |I|, 1\le x_1x_2\le \kk{t}\right\}.
\end{equation}
By the definition of $\F{k}{I}$,
\begin{equation*}
\label{dm2}
    w\in \F{k}{I}\iff (\av{w}{J}, \av{w^{-1}}{J},|J|)\in\omega\text{\space for any subinterval \space} J\subset I.
\end{equation*}
We will also need planar projections of $\omega$ for fixed $t\in(0,|I|]$: 
    \begin{equation}
    \label{dm3}
        \omega^t=\{(x_1,x_2)\in\mathbb{R}^2\colon 1\le x_1x_2\le k(t)\}.
    \end{equation}
We now introduce a family of curves that will play a central role in the proof of the main theorem. 
    For any fixed $t>0$ define
    \begin{equation*}
    \label{c2}
        \Cn{\tau}{t}{}=(\Cn{\tau}{t}{1}{},\Cn{\tau}{t}{2}{}=\left(1+\frac{\tau}{t}(\rr{\tau}-1), 1+\frac{\tau}{t}\left(\kr{\tau}-1\right)\right),\quad\tau\in[0,|I|],
    \end{equation*}
and let $\gamma^t$ be the restriction of the curve $\Gamma^t$ to the interval $[0,t]$:
    \begin{equation*}
    \label{c3}
        \cn{\tau}{t}{}=\Cn{\tau}{t}{},~\tau\in[0,t].
    \end{equation*}
Let us collect key properties of the curves $\Gamma^t$ and $\gamma^t$ that will be needed in the remainder of the proof.
    \begin{lemma}
    \label{cL1}
        The functions $\Cnn_1^t$ and $\Cnn_2^t$ satisfy the relations 
        \begin{equation}
        \label{c31}
            (\Cnn_1^t)'>0,
        \end{equation}
        \begin{equation}
        \label{c32}
            (\Cnn_2^t)'<0,
        \end{equation}
        and
        \begin{equation}
        \label{c33}
            \bigg[\frac{\Cnn_2^t-1}{\Cnn_1^t-1}\bigg]'>0.
        \end{equation}
    \end{lemma}
    \begin{proof}
        It is enough to consider the case $t=1.$ Write $a=\Gamma^1_1-1,$ $b=\Gamma^1_2-1.$ Using\eqref{r2} we get
\begin{align*}
\label{c4}
A'(\tau)&=\de{\tau^2 (k(\tau)-1)}=\tau\Big(\kr{\tau}-1\Big)+\tau\left(\rr{\tau}-1\right)=a+b.
\end{align*}
Thus, $a'+b'=A''$ and from~\eqref{c3}, $a'+\frac{r^2}k\,b'=0.$ Therefore,
\eq[un]{
a'=\frac{A''r^2}{r^2-k},\qquad b'=-\frac{A''k}{r^2-k}.
}
Since $A''>0$ and, immediately from~\eqref{r2}, $r^2>k,$ we have established~\eqref{c31} and~\eqref{c32}. 

To show~\eqref{c33}, we compute: Now, \eqref{r2} and \eqref{r3}
\begin{align*}
\Big(\frac ba\Big)'=\frac{b'a-a'b}{a^2}=\frac{b'}{a^2}\,\Big(a+\frac{r^2}k\,b\Big)=\frac{b'\tau}{a^2}\,\Big(2r-1-\frac{r^2}k\Big).
\end{align*}
Since $b'<0$ and, from~\eqref{r2}, $2r-1-\frac{r^2}k=-\tau\frac{k'}k<0,$ this completes the proof. 
    \end{proof}
\begin{figure}
    \centering
    \includegraphics[height=7cm]{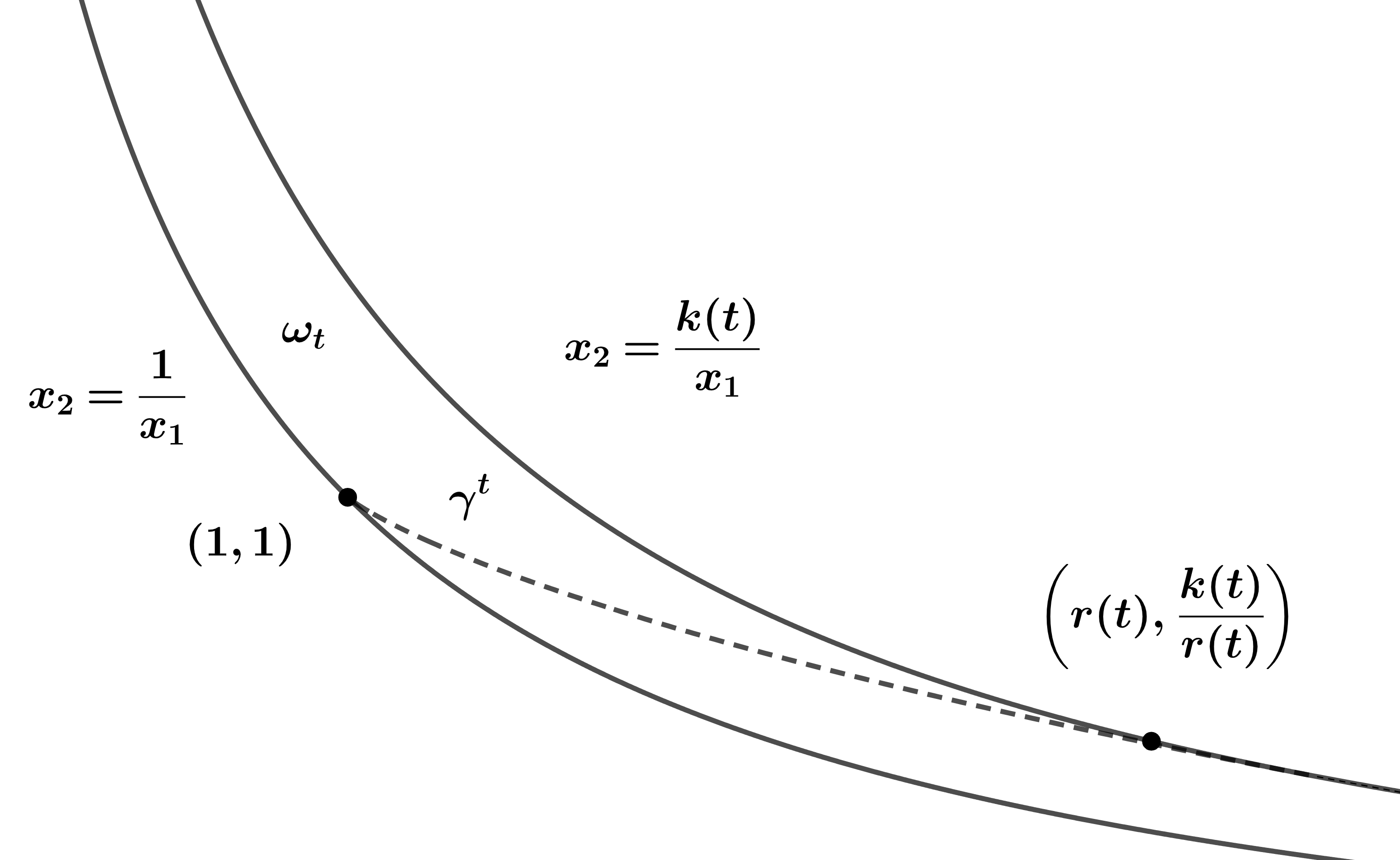}
    \caption{The curve $\gamma_t$}
    \label{figure 1}
\end{figure}

\begin{lemma}
\label{cL2}
    If $0<t_1<t_2$, then the curve $\Cnn^{t_1}$ lies below the curve $\Cnn^{t_2}$
\end{lemma}
\begin{proof}
Observe that 
    $$
    \frac{\Cn{\tau}{t_2}{2}-1}{\Cn{\tau}{t_2}{1}-1}=\frac{\Cn{\tau}{t_1}{2}-1}{\Cn{\tau}{t_2}{1}-1}.
    $$
    Therefore, for $\tau>0$ the point $\Cnn^{t_2}(\tau)$ lies on the line segment connecting the points $(1,1)$ and $\Cnn^{t_1}(\tau).$ 
   On the other hand, from inequalities \eqref{c31} and \eqref{c33} we conclude that this line segment lies above the curve $\{\Gamma^{t_1}(s), 0<s<\tau\},$ which completes the proof.
\end{proof}
\begin{lemma}
\label{cL3}
    For any $t\in(0,|I|)$ the function
    \begin{equation*}
    \label{cL3e1}
        \tau\mapsto\Cn{\tau}{t}{1}\Cn{\tau}{t}{2},\quad 0\le \tau\le |I|,
    \end{equation*}
    is strictly increasing on 
    $[0,t]$ 
    and strictly decreasing on 
    $[t,|I|]$. 
    Its maximum is attained at 
    $\tau=t$ 
    and equals 
    $\kk{t}$. Consequently, for any $s,t>0$,
    \begin{equation}
     \label{cL3e2}
            \Big(1+\frac{s}{t}\left(\rr{s}-1\right)\Big)\left(1+\frac{s}{t}\left(\kr{s}-1\right)\right)\le\kk{t}.
    \end{equation}
    
\end{lemma}
\begin{proof}
   Write $a=\Gamma_1^{t_1}-1=\frac{\tau}t\,(r-1),$ $b=\Gamma_2^{t}-1=\frac{\tau}t\,(\frac kr-1).$ Just as in~\eqref{un} in the proof of~Lemma~\ref{cL1}, $a'=-\frac{r^2}k\,b'$. Therefore,
   $$
   (\Gamma_1^t\Gamma_2^t)'=b'\Big(-\frac{r^2}k\,(b+1)+a+1\Big)=-b'\Big(\frac{r^2}k-1\Big)\Big(1-\frac \tau t\Big).
   $$
 Since, $b'<0$ and $r^2>k,$ the function $\tau\mapsto \Gamma_1^t(\tau)\Gamma_2^t(\tau)$ has a strict maximum at $\tau=t.$
   
To prove~\eqref{cL3e2}, note that it is simply the statement
$$\Cn{s}{t}{1}\Cn{s}{t}{2}\le\kk{t},$$
which is true since $\Cn{s}{t}{1}\Cn{s}{t}{2}<\Cn{t}{t}{1}\Cn{t}{t}{2}=k(t)$ by the first part of the claim.
\end{proof}
\begin{figure}
    \centering
    \includegraphics[height=7cm]{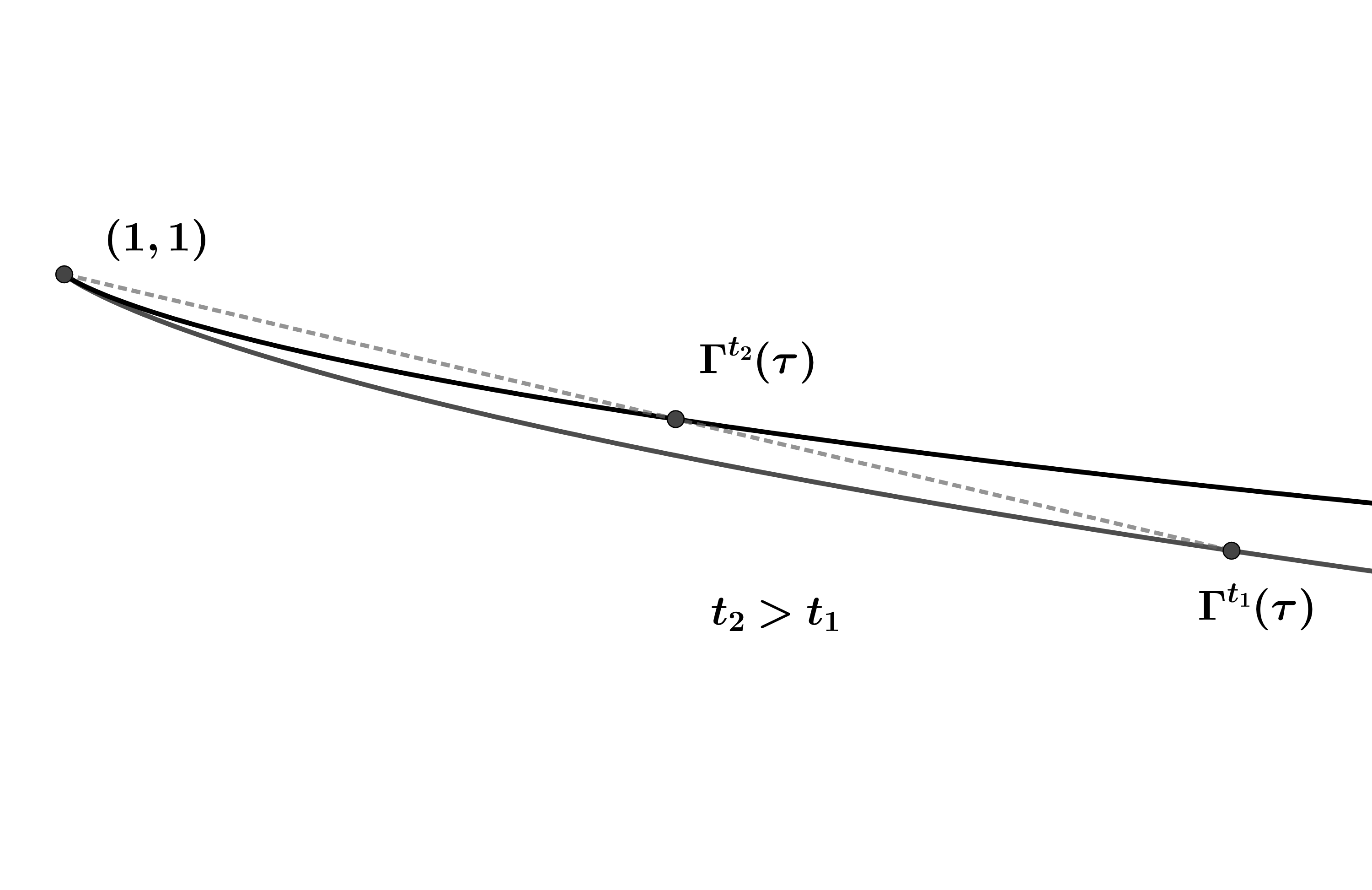}
    \caption{$\Cnn^{t_1}$ vs $\Cnn^{t_2}$}
    \label{figure 2}
\end{figure}
This gives the following corollary.
\begin{corollary}
\label{cC1}
    For any $t>0$ the curve $\Cnn^t$ lies below the hyperbola $x_1x_2=\kk{t}$; the two are tangent when $x_1=\rr{t}$, which is their only point in common. Consequently, for each $C\in[1,\kk{t}]$ the curve $\cnn^t$ intersects the hyperbola $x_1x_2=C$ exactly once.
\end{corollary}
\subsection{Special functions $U$ and $\mathcal{T}$}
\begin{lemma}
\label{cL4}
    For any $t>0$ and point $(y_1,y_2)\in \omega^t$ there exist unique $\tau=\TT{y}{t}\in[0,t]$ and $u=\UU{y}{t}\in(0,\infty)$ such that 
    \begin{equation*}
    \label{cL4e1}
       y_1=u\cn{\tau}{t}{1}
    \end{equation*}
   and
    \begin{equation*}
    \label{cL4e2}
        y_2=\frac{\cn{\tau}{t}{2}}{u}.
    \end{equation*}
    The functions $U$ and $\mathcal{T}$ satisfy the following homogeneity relations: for any $c>0$
    \begin{equation*}
    \label{cL4e3}
        \mathcal{T}\left(cy_1,\frac{y_2}{c},t\right)=\TT{y}{t}.
    \end{equation*}
    \begin{equation*}
    \label{cL4e4}
        U\left(cy_1,\frac{y_2}{c},t\right)=c\,\UU{y}{t}.
    \end{equation*}
\end{lemma}
\begin{proof}
\begin{figure}
    \centering
    \includegraphics[height=8cm]{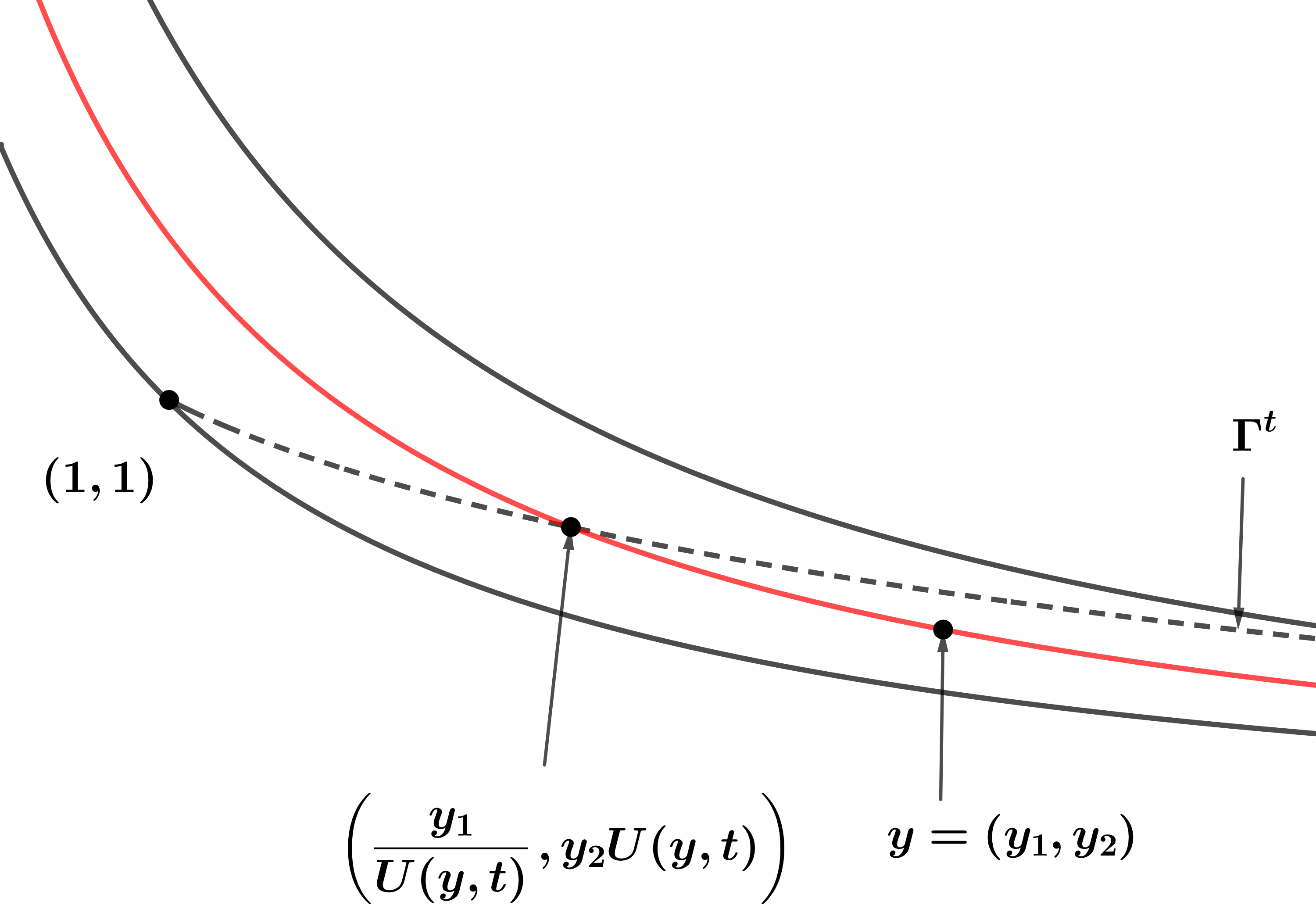}
    \caption{Definition of $U$ and $\mathcal{T}$}
    \label{figure 3}
\end{figure}
Let $C=y_1y_2$. Then the hyperbola $x_1x_2=C$ passes through the point $(y_1,y_2)$ and, by Corollary~\ref{cC1}, intersects the curve $\cnn^t$ at a single point. Call the $x_1$-coordinate of that point $\cnn^t(\mathcal{T})$ and let $U=\frac{y_1}{\cnn^t_1(\mathcal{T})}$.
Therefore,
$$
y_1y_2=\cn{\mathcal{T}}{t}{1}\cn{\mathcal{T}}{t}{2}
$$
and
$$
y_2=\frac{\cn{\mathcal{T}}{t}{1}\cn{\mathcal{T}}{t}{2}}{y_1}=\frac{\cn{\tau}{t}{2}}{U},
$$
as claimed.
To show the homogeneity relations note that $(cy_1)(\frac{y_2}{cy_1})=y_1y_2,$  so the point $(cy_1,\frac{y_2}c)$ lies on the same hyperbola as the point $(y_1,y_2).$
\end{proof}
We will proceed by stating the following corollary, which will be useful later. The first two statements follow directly from the definitions of $U$ and $\mathcal{T}.$ The last statement follows from Lemma \ref{cL2}. 
\begin{corollary}
\label{uC1}
    The functions $U$ and $\mathcal{T}$ satisfy the following properties:
    \begin{itemize}
    \item $U(y_1,y_2,t)\in\big[\frac{y_1}{r(t)},y_1\big];$
    \item If $y_1y_2=1,$ then $\mathcal{T}(y_1,y_2,t)=0$ and $U(y_1,y_2,t)=y_1;$
    \item For $0<t_1<t_2$ and any $y\in\omega^t$ such that $y_1y_2>1$ we have
$$
        U(y_1,y_2,t_1)< U(y_1,y_2,t_2).
$$
\end{itemize}
\end{corollary}
\subsection{The main estimate on $\inf w$}
Given $a>0$, define the following restriction of the domain \eqref{dm3}. 
\begin{equation*}
\label{u1}
    \omega_a^t=\{x\in\omega^t:U(x,t)\ge a\}.
\end{equation*}
Our goal is to establish that for any $w\in\F{k}{I}$, $(\av{w}{I},\av{w^{-1}}{I})$ lies in $\omega_{\inf_{I}w}^{|I|}$. To that end, we first prove the following key lemma.  
\begin{lemma}
\label{uL1}
    If $w\in F_2^{k}(I)$ and $\av{w}{I}\av{w^{-1}}{I}<k(|I|),$ then there exists a subinterval $J\subset I$ such that $0<|J|<|I|$ and
    \begin{equation*}
    \label{uLe1}
        U(\av{w}{J},\av{w^{-1}}{J},|J|)\le U(\av{w}{I},\av{w^{-1}}{I},|I|).
    \end{equation*}
\end{lemma}
\begin{figure}
    \centering
    \includegraphics[height=9cm]{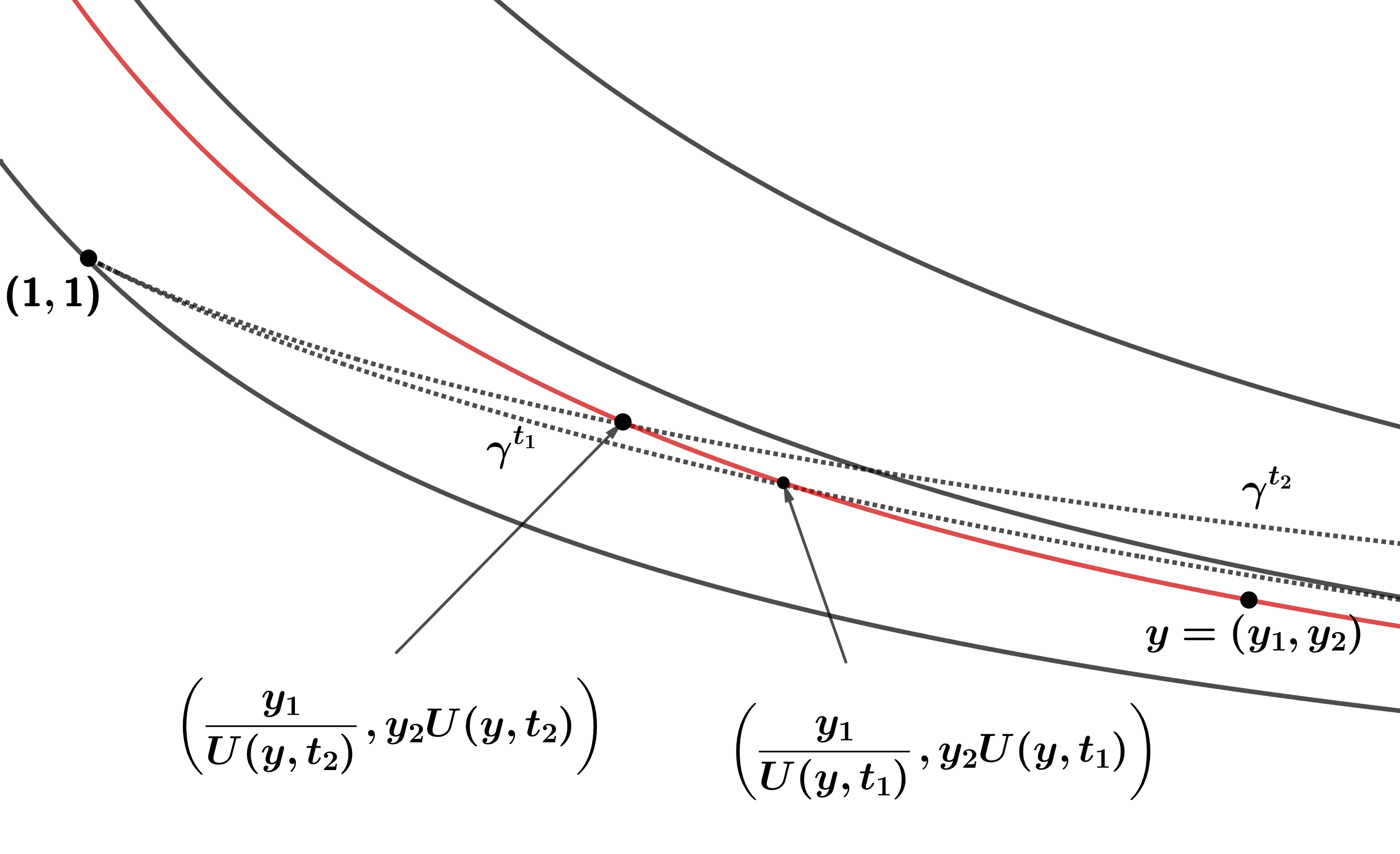}
    \caption{Definition of $U$ and $\mathcal{T}$}
    \label{figure 4}
\end{figure}
\begin{proof}
Let $x_1=\av{w}I,$ $x_2=\av{w^{-1}}I,$ and $|I|=t$. Without loss of generality, we may choose $I=[0,t]$; by multiplying $w$ by a constant if necessary, we may assume $\UU{x}{t}=1$. Then, from Lemma~\ref{cL4},
\begin{equation}
\label{uLe2}
    x_1=\Cn{\tau}{t}{1},~x_2=\Cn{\tau}{t}{2}, \text{~where~} \tau=\TT{x}{t}.
\end{equation} 
If $x_1x_2=1$ then $w=1$ a.e.  on $I.$ Thus, for any $J\subset I$ we have $U(\av{w}J,\av{w^{-1}}J,|J|)=1$.
From now on, assume $x_1x_2>1$. Then $\tau>0$. Since $x_1x_2<\kk{t},$ we have $\tau< t$ and $x_1< \rr{t}$. By continuity of $r$, there exists a $t_0\in(\tau,t)$ such that $tx_1<t_0\rr{t_0}$.
If $\av{w}J\le1$ for some subinterval $J\subset I$ with $0<|J|<|I|,$ then $U(\av{w}J,\av{w^{-1}}J,|J|)\le\av{w}J\le 1$ and~\eqref{uLe2} holds. Thus, let us assume that $\av{w}J>1$ for every $J\subset I$, which means that $w\ge 1$ a.e. on $I$.

Consider the following ratio: 
\begin{equation*}
\label{uLe3}
    \M(J)=\frac{\av{w^{-1}}J-1}{\av{w}J-1},\text{\space\space\space}J\subset I.
\end{equation*}
Let $s\in (0,t)$ and consider the intervals $J_-(s):=(0,s)$ and $J_+(s):=(s,t)$. Then $\M(I)$ lies between $\M(J_-(s))$ and $\M(J_+(s))$ and $\M(I)=\M(J_-(s))$ if and only if $\M(I)=\M(J_+(s))$.
Suppose that for some $s\in(0,t)$ we have 
\begin{equation*}
\label{uLe4}
    \M(J_-(s))=\M(I)=\M(J_+(s)).
\end{equation*}
Set $J=J_-$ or $J=J_+$ in such a way that 
$\av{w}J\leq x_1$. 
Put $\lambda=\frac{x_1-1}{\av{w}J-1},$ then $\lambda\ge1$. Observe that
\begin{equation*}
    \label{uLe5}
    \av{w}J-1=
\frac{x_1-1}{\lambda}=
\frac{\Cn{\tau}{t}{1}-1}{\lambda}=
\frac{\tau}{\lambda t}(\Cn{\tau}{1}{1}-1)).
\end{equation*}
Therefore, $\av{w}J=\Gamma_1^{\lambda t}(\tau)$.
Similarly, $\av{w^{-1}}J=\Gamma_2^{\lambda t}(\tau)$. 
Since $\lambda\ge1$, $\tau\in[0,\lambda t]$. Therefore,
$\mathcal{T}(\av{w}J,\av{w^{-1}}J,\lambda t)=\tau$. Hence, by the third statement in Corollary~\ref{uC1},
\begin{equation*}
\label{uLe6}
    1=U(\av{w}J,\av{w^{-1}}J,\lambda t)\ge U(\av{w}J,\av{w^{-1}}J,|J|),
\end{equation*}
so~\eqref{uLe2} holds. Now assume that the statement $\M(J_-(s))=\M(I)=\M(J_+(s))$ is false for any $s\in(0,t).$ By continuity of $\M(J_{\pm})$ with respect to $s$, we have either
\begin{equation*}
\label{uLe7}
 \M(J_-(s))<\M(I)<\M(J_+(s))\text{\space\space or\space\space}\M(J_+(s))<\M(I)<\M(J_-(s))
\end{equation*}
for all $s\in(0,t)$ simultaneously.
Choosing $s=t_0$ for one bound and $s=t-t_0$ for the other, we conclude that $\M(I)$ lies between $\M(J_-(t_0))$ and $\M(J_+(t-t_0))$.
Since the function $s\mapsto \M(s,s+t_0)$ is continuous on $[0,t-t_0],$ there exists $s\in [0,t-t_0]$
such that for $J=(s,s+t_0)$ we have $\M(J)=\M(I)$. Since $w\ge 1,$ we have
\begin{equation*}
\label{uLe8}
    \av{w}J-1\le\frac{|I|}{|J|}\,(\av{w}I-1)=\frac{t}{t_0}\,(x_1-1)<r(t_0)-1.
\end{equation*}
Let  
$\mu=\frac{\Cn{\tau}{1}{1}-1}{t_0(\av{w}J-1)}=\frac{t(x_1-1)}{t_0(\av{w}J-1)}.$ Then $\mu \ge1$ and
\begin{equation*}
\label{uLe9}
   \av{w}J=1+\frac{\Cn{\tau}{}{1}-1}{\mu t_0}=\Gamma_1^{\mu t_0}.
\end{equation*}
Furthermore,
\begin{align*}
\label{uLe10}
   \av{w^{-1}}J&=1+\M(J)(\av{w}J-1)=1+\M(I)(\av{w}J-1)\\&=1+\bigg[\frac{\av{w^{-1}}I-1}{\av{w}I-1}\bigg](\av{w}J-1)=1+\bigg[\frac{\av{w^{-1}}I-1}{\av{w}I-1}\bigg]\bigg[\frac{\Cn{\tau}{1}{1}-1}{\mu t_0}\bigg]\\&
   =1+\frac{\Cn{\tau}{1}{2}-1}{\mu t_0}=1+\Cn{\tau}{\mu t_0}{2}.
\end{align*}
Since $\tau\in[0,\mu t_0]$, we have $\mathcal{T}(\av{w}J,\av{w^{-1}}J,\mu t_0)=\tau$. Hence, by Corollary~\ref{uC1},
\begin{equation*}
\label{uLe12}
    1=U(\av{w}J,\av{w^{-1}}J,\lambda t_0)\ge U(\av{w}J,\av{w^{-1}}J,|J|).
\end{equation*}
This completes the proof.
\end{proof}
\begin{theorem}
\label{uT1}
    If $w\in \F{k}{I}$, then 
    \begin{equation}
    \label{uT1e1}
        \inf_I w\le U(\av{w}I,\av{w^{-1}}I,|I|).
    \end{equation}
    Equivalently, the point $(\av{w}I,\av{w^{-1}}I )$ lies in $\omega_{\inf_I w}^{|I|}$.
\end{theorem}
\begin{proof}
  For any subinterval $J\subset I,$ let
  \begin{equation*}
      \label{uT1e2}
      L(J)=U(\av{w}J,\av{w^{-1}}J,|J|).
  \end{equation*}
Assume $w\in \F{k}{I}$. It is enough to consider the case $L(I)=1$. Let
\begin{equation*}
\label{uT1e3}
    L=\inf\{L(J)\colon J\subset I\}.
\end{equation*}
Clearly $L\le L(I)=1$. Pick a sequence of intervals $J_n$ such that $L(J_n)\rightarrow L$ and $D:=\lim |J_n|$ is minimal.
To get a contradiction, assume $D>0$.
Then there exists an interval $J\subset I$ such that $|J|=D$, and a subsequence $\{J_{n_{k}}\}_k$ that converges to $J$. So
\begin{equation*}
    \label{uT1e4}
    L(J)=\lim L(J_{n_{k}})=L
\end{equation*}
by the continuity of the function $U$. By Lemma ~\ref{uL1}, there exists an interval $\tilde{J}\subset J$ such that $L(\tilde{J})\le L(J)$ and $0<|\tilde{J}|<|J|=D,$ contradicting the minimality of $J.$  Therefore, $D=0$. 
Hence for any $\delta>0$ there exists an interval $J\subset I$ such that $L(J)\le 1+\delta$ and $r(|J|)\le1+\delta$ (here we used the continuity of $r$ at 0, Lemma~\ref{rL2}).
Notice that 
\begin{equation*}
\label{uT1e5}
  \inf_{I}w\le\inf_{J}w\le\av{w}J\le L(J)r(|J|)\le (1+\delta)^2,   
\end{equation*}
where we used the first statement in Corollary~\ref{uC1}. Taking the limit $\delta \rightarrow 0$ gives $\inf_{I}w\le 1$.
\end{proof}
\subsection{The cutout function and the final part of the proof}
The following definition is from~\cite{sz1}.
\begin{definition}
  Let $I=[t_1,t_2]\subset\mathbb{R}$ be an interval and let $E\subset\mathbb{R}$ be a set of positive measure such that $d:=|I\setminus E|>0$. Define a function $h_{I,E}:I\to\mathbb{R}$ by
  $$h_{I,E}(s)=|[t_1,s]\setminus E|.$$
  This function is increasing and Lipchitz. It has the image $[0,d].$ 
The pre-image  $h^{-1}_{I,E}(\tau)$ is a one-point set for all but at most countably many $\tau\in[0,d].$

For a weight $w\colon I\to \mathbb{R}$ let
$$
  \col_{I,E}w(\tau)=w(h^{-1}_{I,E}(\tau)),\text{\space}\tau\in[0,d].
$$
The function $\col_{I,E}w$ is called the cutout of $w$ along $E.$ This function is well defined for all $\tau\in [0,d]$, except possibly on a countable set, and is measurable.
\end{definition}
\begin{lemma}
\label{uL2}
   Fix $a\in(0,\infty)$. Take $x\in\omega_a^t$ and $s\in(0,t)$, and let 
   \begin{equation*}
   \label{uL2e1}
       y=(y_1,y_2)=(a,a^{-1})+\frac{t}{s}(x-(a,a^{-1})).
   \end{equation*}
 Then $y_1y_2\leq k(s)$.
\end{lemma}
\begin{proof}
    Without loss of generality, assume $a=1$.
Then 
\begin{equation*}
\label{uL2e2}
    y=\left(1+\frac{t}{s}(x_1-1),1+\frac{t}{s}(x_2-1)\right).
\end{equation*}
We proceed with considering two cases.
First, assume $x_1\in[1,\rr{t}]$. Then the point $x$ lies below the curve $\cnn^t$. This means that there exists $\tau\in[0,t]$ such that 
\begin{equation*}
\label{uL2e3}
    x_1=\Cn{\tau}{t}{1},\quad x_2\leq\Cn{\tau}{t}{2}.
\end{equation*}
Then 
\begin{align*}
    y_1&=1+\frac{t}{s}(x_1-1)=1+\frac{t}{s}\left(\frac{\tau}{t}\left(\rr{\tau}-1\right)\right)=\Cn{\tau}{s}{1}
\end{align*}
and
\begin{align*}
     \label{uL2e5}
    y_2&\le 1+\frac{t}{s}(x_2-1)\le1+\frac{t}{s}\left(\frac{\tau}{t}\left(\frac{k(\tau)}{\rr{\tau}}-1\right)\right)=\Cn{\tau}{s}{2}=\frac{k(s)}{\Cn{\tau}{s}{1}}=\frac{k(s)}{y_1},
\end{align*}
so the claim is proved in this case. Assume now that $x_1>\rr{t}$. 
We have 
$$
\label{uL2e7}
    y_2=1+\frac{t}{s}(x_2-1)\le 1+\frac{t}{s}\left(\frac{\kk{t}}{x_1}-1\right),
$$
thus, the claim will be proved if we show that 
\begin{equation*}
\label{uL2e8}
    \Big[1+\frac{t}{s}\Big(\frac{\kk{t}}{x_1}-1\Big)\Big]\Big[1+\frac{t}{s}(x_1-1)\Big]\le k(s).
\end{equation*}
The derivative of the left-hand side of this expression with respect to $x_1$ is 
$
\frac ts\,\Big(1-\frac ts\Big)\Big(1-\frac{k(t)}{x_1^2}\Big).
$
Since $t>s,$ $1-\frac ts<0.$ Since $x_1>r(t),$ we have $1-\frac{k(t)}{x_1^2}>1-\frac{k(t)}{r^2(t)}>0.$ Thus, the derivative is negative and it suffices to prove that
$$
\Big[1+\frac{t}{s}\Big(\frac{\kk{t}}{r(t)}-1\Big)\Big]\Big[1+\frac{t}{s}(r(t)-1)\Big]\le k(s),
$$
which is inequality~\eqref{cL3e2} of Lemma~\ref{cL3} with $s$ and $t$ interchanged. The proof is complete.
\end{proof}
\begin{prop}
\label{uP1}
Take $w\in  \F{k}{I}$ and assume that $\inf_I w\ge a$ for some $a\in(0,\infty)$.
Let $E\subset I$ be such that $w=a$ a.e. on $E$. Let $d=|I\setminus E|$. Let $v=\col_{I,E}w.$ Then $v\in F_2^k([0,d]).$
\end{prop}
\begin{proof}
 Let $\tau_1,\tau_2\in[0,d],\tau_1<\tau_2$, $|J|=[\tau_1,$ $\tau_2],$ and $s=|J|.$ 
 Without loss of generality assume $h_{I,E}^{-1}(\tau_1)$ and  $h_{I,E}^{-1}(\tau_2)$ are one-point sets. Let $\tilde{J}=h_{I,E}^{-1}(J)$, $t=|\tilde{J}|$, and 
 \begin{equation*}
 \label{uL2e12}
     x=(x_1,x_2)=(\av{w}{\tilde{J}},\av{w^{-1}}{\tilde{J}}).
 \end{equation*}
From Theorem \ref{uT1}, $U(x,t)\ge \inf_{\tilde{J}}w\ge a,$ meaning $x\in\omega_a^t.$ Let
 \begin{equation*}
     \label{uL2e13}
     y=(y_1,y_2)=(\av{v}J,\av{v^{-1}}J).
 \end{equation*}
Observe that
\begin{align*}
\label{uL2e14}
    \av{v}J=\frac{1}{|J|}\left(\int_J a+\int_{J}(v-a)\right)=\frac{1}{|J|}\left(\int_J a+\int_{\tilde{J}}(w-a)\right)=a+\frac{t}{s}(\av{w}{\tilde{J}}-a).
\end{align*}
 Similarly,
 \begin{equation*}
     \label{uL2e15}
     \av{v^{-1}}{J}=a^{-1}+\frac{t}{s}\big(\av{w^{-1}}{\tilde{J}}-a^{-1}\big).
 \end{equation*}
 Thus,
 \begin{equation*}
 \label{uL2e16}
     y=(a,a^{-1})+\frac{t}{s}(x-(a,a^{-1})).
 \end{equation*}
By Lemma \ref{uL2},
 $y_1y_2\leq k(s)$.
 Therefore, $y\in\omega^s$, which completes the proof.
\end{proof}
The symmetric statement is also true.
\begin{prop}
\label{uP2}
Take $w\in  \F{k}{I}$ and assume that $\sup w\ge b$ for some $b\in(0,\infty)$.
Let $E\subset I$ be such that $w=b$ a.e. on $E$. Let $d=|I\setminus E|$. Let $v=\Cut_{I,E}w.$ Then $v\in F_2^k([0,d]).$
\end{prop}
Propositions \ref{uP1} and \ref{uP2} together give the following combined statement.
\begin{prop}
\label{uP3}
Take $w\in  \F{k}{I}$ and assume that $a\le w\le b$ for some $a,b\in(0,\infty)$.
Let $E_a\subset I$ and $E_b\subset I$  be such that $w=a$ a.e. on $E_a$ and $w=b$ a.e on $E_b$. Let $E=E_a\cup E_b$ and $d=|I\setminus E|$. Let $v=\Cut_{I,E}w.$ Then $v\in F_2^k([0,d]).$ 
\end{prop}
We need one last lemma to finish the proof of Theorem~\ref{MT1} in the smooth convex case. 
\begin{remark}
The lemma is a simple statement that the truncation of a weight does not increase its $A_2$-oscillation on any interval. The proof is stated for an interval, but works for any set in any dimension. The result itself is not new -- for instance, there is an elementary proof in \cite{rvv}, or it can be derived in one line from Lemma~4.2 from~\cite{slavin}. However, the argument we give is different (the idea is due to Pavel Zatitskii) and is similar in spirit to the one-line proof that the $\BMO$-oscillation does not increase under truncation, which is obvious once one writes that oscillation in this manner:
\eq[14]{
\av{\varphi^2}J-\av{\varphi}J^2=\frac12\,\frac{1}{|J|^2}\int_J\int_J(\varphi(x)-\varphi(y)^2\,dx\,dy.
}
The direct analog of~\eqref{14} for $A_2$ is this:
\begin{equation}
 \label{uL3e4}
   \av{w}{J}\av{w^{-1}}{J}=\frac12\,\frac{1}{|J|^2}\int_J\int_J\Big(\frac{w(x)}{w(y)}+\frac{w(y)}{w(x)}\Big)dxdy.
\end{equation}
\end{remark}

    \begin{lemma}
    \label{uL3}
 Let $w$ be a weight on $J$ such that both $w$ and $w^{-1}$ are integrable. Take $c,d\in\mathbb{R}$ such that $0\le c<d$. Let $w_{c,d}$ be the truncation of $w$ at levels $c$ and $d$:
 \begin{equation}
  \label{uL3e1}
    w_{c,d}(s) = 
    \begin{cases}
  c, &\text{\space if  } w(s)<c;\\
  w(s), &\text{\space if  } c\le w(s)\le d;\\
 d, &\text{\space if  } w(s)> d.
\end{cases} 
 \end{equation}
Then
 \begin{equation}
  \label{uL3e2}
     \av{w_{c,d}}{J}\av{w^{-1}_{c,d}}{J}\le \av{w}{J}\av{w^{-1}}{J}.
 \end{equation}

\end{lemma}
\begin{proof}
We will show that the one-sided truncation from above at level $d$ does not increase the $A_2$-oscillation, i.e., prove~\eqref{uL3e2} for $c=0.$ The two-sided truncation is just the lower truncation applied to the upper truncation, and the argument for the lower truncation is virtually identical. We will also assume, without loss of generality, that $w>0$ on $J.$

Let
 \begin{equation*}
  \label{uL3e3}
     E_d=\{s\in J: w(s)> d\}.
 \end{equation*}
 If $x,y\in E_d$, then
\begin{equation*}
 \label{uL3e5}
    \frac{w_{0,d}(x)}{w_{0,d}(y)}+\frac{w_{0,d}(y)}{w_{0,d}(x)}=2\le \frac{w(x)}{w(y)}+\frac{w(y)}{w(x)}.
\end{equation*}
If $x,y\notin E_d $, then
\begin{equation*}
 \label{uL3e6}
    \frac{w_{0,d}(x)}{w_{0,d}(y)}+\frac{w_{0,d}(y)}{w_{0,d}(x)}=\frac{w(x)}{w(y)}+\frac{w(y)}{w(x)}.
\end{equation*}
Finally, if $x\in E_d$ and $y\notin E_d$, then $w(y)\le d\le w(x)$ and
\begin{equation*}
 \label{uL3e7}
   \frac{w_{0,d}(x)}{w_{0,d}(y)}+\frac{w_{0,d}(y)}{w_{0,d}(x)}=\frac{d}{w(y)}+\frac{w(y)}{d}\le\frac{w(x)}{w(y)}+\frac{w(y)}{w(x)},
\end{equation*}
where we used the fact that the expression in the middle is increasing in $d$ when $d\ge w(y).$
In all three cases, integrating both sides and using~\eqref{uL3e4} we obtain~\eqref{uL3e2} for $c=0.$
\end{proof}
We can now finish the proof of Theorem \ref{MT1} in the smooth convex case.
\begin{proof}
   Without loss of generality, assume
$I=[0,t].$ Take any interval $J=[t_1,t_2]\subset I.$ We need to show that
\eq[716]{
\av{w^*}J\av{(w^*)^{-1}}J\le k(t_2-t_1).
}
Let $c=w^*(t_2)$ and $d=w^*(t_1);$ thus $c\le d.$ Let us take the truncation $w_{c,d}$ in place of $w.$ This truncation does not affect~\eqref{716} and by Lemma~\ref{uL3}, $w_{c,d}\in F_2^k(I).$

Pick sets $E_c\subset I, E_d \subset I$ such that $w=c$ on $E_c$, $w=d$ on $E_d$, $|E_c|=t-t_2,$ and $E_d=t_1$. Let $E=E_c\cup E_d$ and $\delta=|I\setminus E|$. From Proposition \ref{uP3}, $v=\Cut_{I,E}w$ lies in $F_2^k([0,\delta])$. Since the distribution of $w^*$ on $[t_1,t_2]$ coincides with the distribution of $v$ on $[0,\delta],$ we get~\eqref{716}.
\end{proof}

\section{The proof of the main theorem in the general convex case}
\label{gen1}
In this section we will prove Theorem \ref{MT1} in the general convex case. Specifically, we assume that $k$ is an increasing continuous function on $[0,|I|]$ and $\AAA(t)=t^2(k(t)-1)$ is convex on $[0,|I|],$ with no assumptions on the differentiability of $k.$

We will use the following mollification procedure, which is a minor modification of the one in~\cite{sz1}.
\label{GenSet}
Fix an interval $I.$ Let $T=|I|.$ Since $A$ is defined on the finite interval $[0,T]$, while our procedure works for the half-line, we first extend $A$ to $[0,\infty).$ Specifically, note that $A$ is convex and continuous on $[0,T]$. Therefore, it has the left derivative at $t=T;$ denote it by $A'(T)$. Let 
\eq[new_a]{
\hat A(t)=
\begin{cases}
A(t),& 0\le t\le T;\\
A'(T)(t-T)+A(T),& T<t<\infty.
\end{cases}
} 
Clearly, $\hat{A}$ is continuous and convex on $[0,\infty).$ This also extends $k$ to the half-line:
\eq[new_k]{
\hat k(t)=
\begin{cases}
k(t),& 0\le t\le T;\\
\frac1{t^2}\,\hat{A}(t)+1,& T<t<\infty.
\end{cases}
}
Here is the procedure from~\cite{sz1}. Pick a function $\Psi\in C_{0}^{\infty}[0,\infty)$ such that $\Psi\ge 0$, $\int_0^\infty u\Psi(u)du=1$, and $\int_0^\infty\Psi(u)du=1.$
For a continuous function $F$ defined on $[0,\infty),$ consider the following multiplicative convolution:
\begin{equation}
\label{Ge1}
    F\star\Psi(t)=\int_0^{\infty}\frac{F(v)}{t}\,\Psi\Big(\frac{v}{t}\Big)\,dv, \quad 0<t< \infty.
\end{equation}
Changing the variable we have
\eq[618]{
F\star\Psi(t)=\int_0^\infty F(tu)\Psi(u)\,du.
}
From~\eqref{Ge1} we see that $F\star\Psi$ is smooth on $(0,\infty).$ Observe that~\eqref{618} coupled with Jensen's inequality implies that if $F$ is convex on $(0,\infty),$ then so is $F\star \Psi$: 
\eq[619]{
F\star\Psi(t)=\int_0^\infty F(tu)\Psi(u)\,du\ge F\left(\int_0^\infty tu \Psi(u)\,du\right)=F(t).
}
Take a sequence of such functions $\{\Psi_n\}$ whose supports converge to the set $\{1\}.$ Then $F\star \Psi_n\to F$ pointwise on $(0,\infty).$

Now, define a sequence of moduli $\{k_n\}$ on $[0,\infty)$ as follows:
\begin{equation}
\label{Ge3}
k_n(t)=
\begin{cases}
\frac{1}{t^2}(\hat{A}\star\Psi_n(t)+\frac{t^3}{n})+1,& t>0,\\
1,& t=0.
\end{cases}
\end{equation}
We now state Lemma~3.1 from~\cite{sz1} as it applies to our definition~\eqref{Ge3} (this is done simply by replacing $\xi_n$ in that lemma with $\sqrt{k_n-1}$).
\begin{lemma}[\cite{sz1}]
\label{GL1}
The functions $k_n$ have the following properties:
     \begin{enumerate}
         \item $k_n$ is smooth, $k_n>1$, and $k'_n>0$ on $(0,\infty)$;
         \item $k_n(t)\rightarrow 1$ as $t\rightarrow 0^{+}$;
         \item $A_n\colon t\mapsto t^2(k_n(t)-1)$ is convex on $[0,\infty)$ and $A''_n(t)>0$ for $t>0$;
         \item $k_n\ge \hat{k}$ on $[0,\infty)$;
         \item $k_n\rightarrow \hat{k} $ as $n\to\infty$ pointwise on $[0,\infty).$
     \end{enumerate}
\end{lemma} 

We can now prove Theorem \ref{MT1} in the general convex case.
\begin{proof}
Pick $w\in F_2^{k}(I)$. Since $k_n\ge k,$ we have $w\in F_2^{k_{n}}(I)$. From the smooth convex case, $w^*\in F_2^{k_{n}}(I)$, meaning
$$
\av{w^*}J\av{(w^*)^{-1}}J\le k_n(|J|),
$$
for any subinterval $J$ of $I.$ Taking the limit as $n\rightarrow\infty$, we obtain $w^*\in F_2^{k}(I)$.
\end{proof}

\section{The proof of the main theorem in the general case}
\label{gen2}
Here we eliminate the convexity condition on $\AAA$ and establish Theorem \ref{MT1} in full generality. To do so, we need the following lemma from \cite{sz1}.
\begin{lemma}[\cite{sz1}]
\label{FN}
    Let $\xi$ be a $\VMO$-modulus on $[0,T]$ and take $t_0 \in (0,T)$. Then for any $\delta>0,$ there exists a $\VMO$-modulus $\tilde{\xi}$ on $[0,T]$ such that $\tilde{\xi}\ge \xi$ on $[0,T]$, $\tilde{\xi}(t_0)\le \xi(t_0)+\delta$, and the function $t\mapsto t\tilde{\xi}(t)$ is convex on $[0,T]$. 
\end{lemma} 
Replacing $\xi$ with $\sqrt{k-1}$ and $\tilde{\xi}$ with $\sqrt{\tilde{k}-1}$ gives the following immediate corollary.
\begin{corollary}
\label{cor6}
    Let $k$ be an $F_2$-modulus on $[0,T]$ and take $t_0\in (0,T)$. Then for any $\delta>0$, there exists an $F_2$-modulus $\tilde{k}$ on $[0,T]$ such that $\tilde{k}\ge k$ on $[0,T]$, $\sqrt{\tilde{k}(t)-1}\le \sqrt{k(t_0)-1}+\delta,$ and the function $t\mapsto t\sqrt{\tilde{k}(t)-1}$ is convex on $[0,T].$
\end{corollary}
We are now in a position to prove Theorem \ref{MT1} as stated.
\begin{proof}
    Without loss of generality, let $I=[0,T].$ Our goal is to show that for any subinterval $J\subset I$ we have
$$
\av{w^*}{J}\av{(w^*)^{-1}}{J}\le k(|J|).
$$
If $|J|=|I|=T$, then
$$
\av{w^*}{J}\av{(w^*)^{-1}}{J}=\av{w}{J}\av{w^{-1}}{J}\le k(|J|).
$$
Assume now that $J$ is a proper subset of $I$ and let $t_0=|J|\in(0,T).$ Fix $\delta>0$. Take the modulus $\tilde{k}$ given by  Corollary~\ref{cor6}. Since $\tilde{k}\ge k,$ we have $w\in F_2^{\tilde{k}}(I);$ since the function $t\to t\sqrt{\tilde{k}(t)-1}$ is convex, so is the function $t\to t^2(\tilde{k}(t)-1)$. By the general convex case of Theorem~\ref{MT1} already proved, 
$$
\av{w^{*}}{J}\av{(w^{*})^{-1}}{J}\le \tilde{k}(t_0)\le (\sqrt{k(t_0)-1}+\delta)^2+1= k(t_0)+\delta^2+2\delta\sqrt{k(t_0)-1}.
$$ 
Letting $\delta\to0$ completes the proof. 
\end{proof}

\section{Proofs of Theorem \ref{IT2} and Corollary~\ref{IT1}}
\label{proofs}
In this section, we will provide the proofs of the two elementary results from in Section~\ref{I}.
To prove Theorem~\ref{IT2}, we need the following two results.
\begin{theorem}[\cite{slavin}]
\label{oPT2}
 Let $I$ be an interval. Assume that $w=e^{\varphi}\in A_{\infty}(I)$ and let $C=[e^{\varphi}]_{A_{\infty}(I)}.$ For all $C>1,$
 \begin{equation*}
     \|\varphi\|_{\BMO(I)}\le|\varepsilon^-|,
 \end{equation*}
 where $\varepsilon^-$ is the unique negative solution of the equation $e^{-{\varepsilon}}=C(1-\varepsilon).$
\end{theorem}
\begin{theorem}[\cite{svas}]
    \label{oPT1}
 Let $I$ be an interval and $0<\varepsilon<1$. Suppose $\varphi\in\BMO(I)$ and $\|\varphi\|_{\BMO(I)}\le\varepsilon$. Then
 \begin{equation*}
     \av{e^{\varphi}}{J}\le \frac{e^{-\ve}}{1-\ve}\,e^{\av{\varphi}{J}}\quad\text{for all $J\subset I$}.
 \end{equation*} 
\end{theorem}

\begin{proof}[Proof of Theorem \ref{IT2}]To prove the first statement, 
take $w\in F_2^{k}(I)$ and let $\varphi=\log w$. Observe that for any $J\subset I,$
    $$[e^\varphi]_{A_{\infty}(J)}\le [e^\varphi]_{A_{2}(J)}\le \kk{|J|}.$$
  Let $\nu(|J|)$ be the unique positive solution of  the equation and $e^{\nu(|J|)}=\kk{|J|}(1+\nu(|J|)).$ By Theorem~\ref{oPT2},
    $$
    \|\varphi\|_{\BMO(J)}\le \nu(|J|)=\sqrt{k(|J|)-1}\,\,\left(\frac{\nu^2(|J|)}{k(|J|)-1}\right)^{1/2}.
    $$
   Since $e^p-1-p\ge \frac16\,(p^2+p^3)$ for any $p\ge0,$ we have
   $$
   \frac{\nu^2(|J|)}{k(|J|)-1} =\frac{\nu^2(|J|)}{\frac{e^{\nu(|J|)}}{1+\nu(|J|)} -1}=\frac{\nu^2(|J|)+\nu^3(|J|)}{e^{\nu(|J|)}-1-\nu(|J|)}\le 6,
   $$
Therefore,
   $$
   \|\varphi\|_{\BMO(J)}\le C\sqrt{k(|J|)-1}
   $$
with $C=\sqrt6,$ as claimed.

To prove the second statement, define a function $M$ on $[0,1)$ by
\eq[67]{
M(s)=\frac{e^{-2s}-(1-s)^2}{(1-s)^2s^2}.
}
It is easy to check that $M$ is increasing and has the range $[1,\infty)$.
Since $\|\varphi\|_{\mathcal{L}_{\xi}(I)}=1,$ for each $J\subset I$ we have
$\|\lambda\varphi\|_{\BMO(J)}\le\lambda\xi(|J|)\le \lambda\xi(|I|)<1.$ By Theorem \ref{oPT1},
\begin{align*}
    \av{e^{\lambda \varphi}}{J}\av{e^{-\lambda \varphi}}{J}
    &\le\Big(\frac{e^{-\lambda\xi(|J|)}}{1-\lambda\xi(|J|)}\Big)^2=1+\lambda^2\,\frac{\big(\frac{e^{-\lambda\xi(|J|)}}{1-\lambda\xi(|J|)}\big)^2-1}{\lambda^2\xi^2(|J|)}\,\xi^2(|J|)\\
    &=1+\lambda^2 M\big(\lambda\xi(|J|)\big)\,\xi^2(|J|)\le1+\lambda^2M\big(\lambda\xi(|I|)\big)\,\xi^2(|J|).
\end{align*}
This
completes the proof.
\end{proof}
 
\begin{proof}[Proof of Corollary \ref{IT1}]
  Suppose first that
  \begin{equation*}
      \label{oP1}
      \liminf_{t\to0^+} \frac{\sqrt{\kk{t}-1}}{t}>0.
  \end{equation*}
By Proposition~\ref{oPL2}, the space $\Xicl{\sqrt{k-1}}(I)$ contains non-constant functions. Take one such function $\varphi$ with $\|\varphi\|_{\Xicl{\sqrt{k-1}}(I)}=1.$ Let $\lambda\in\big(0,\frac1{\sqrt{k(|I|)-1}}\big)$ be such that the product $\lambda^2 M\big(\lambda\sqrt{k(|I|)-1}\big)=1,$ where $M$ is given by~\eqref{67}. Let $w=e^{\lambda\varphi}.$ Then $w$ is non-constant and by the second statement in Theorem \ref{IT2},
  $$
  w\in F_2^{1+(\sqrt{k-1})^2}(I)=F_2^k(I).
  $$
Now suppose that
  $$
  \liminf_{t\rightarrow 0^+} \frac{\sqrt{k(t)-1}}{t}=0.
  $$
Take any $w\in F_2^{k}(I).$ By the first statement in Theorem \ref{IT2}, we have $\log w\in \mathcal{L}_{\sqrt{k-1}}(I)$. Furthermore, by Proposition~\ref{oPL2}, $\log(w)$ is constant. Hence, $w$ is constant.
\end{proof}

\end{document}